\documentclass[10pt]{paper}
\usepackage[english]{babel}

\newtheorem{theorem}{Theorem}
\newtheorem{lemma}{Lemma}

\title{Binary and ternary quasi-perfect codes with small dimensions}
\date{}
\begin{document}
\maketitle

\begin{center}
{\large Tsonka Baicheva, Iliya Bouyukliev, Stefan Dodunekov} \\
{\it Institute of Mathematics and Informatics\\
Bulgarian Academy of Sciences, Bulgaria}\\
{\large Veerle Fack}\\
{\it Vakgroep Toegepaste Wiskunde en Informatica, Universiteit
Gent, Belgium}\\
\end{center}

\section{Introduction}

Let $F_q^n$ be the $n$-dimensional vector space over the finite
field with $q$ elements $GF(q)$. A {\it linear code} $C$ is a
$k$-dimensional subspace of $F_q^n$. For $x,y\in F_q^n$ let
$d(x,y)$ denote the Hamming distance between $x$ and $y$, which is
equal to the number of positions where $x$ and $y$ differ. The
minimum Hamming distance for a code $C$ is defined by
$$d(C)=\min_{c_1,c_2\in C,c_1\ne c_2} d(c_1,c_2)$$ and the Hamming
weight $w(x)$ of a vector $x\in F_q^n$ is defined by
$$w(x)=d(x,{\bf 0})$$ where ${\bf 0}$ is the all zero vector. The packing radius
$e(C)$ of the code is $$e(C)=\left\lfloor{{d(C)-1}\over
2}\right\rfloor$$ and this is the maximum weight of successfully
correctable errors. The ball of radius $t$ around a word $y\in
F_q^n$ is defined by $$\{x\vert x\in F_q^n, d(x,y)\leq t\}.$$ Then
$e(C)$ is the largest possible integer number such that the balls
of radius $e(C)$ around the codewords are disjoint. The covering
radius $R(C)$ of a code $C$ is defined as the least possible
integer number such that the balls of radius $R(C)$ around the
codewords cover the whole $F_q^n$, i.e.
$$R(C)=\max_{x\in F_q^n}
\min_{c\in C} d(x,c).$$ With these notations a $q$-ary linear code
of length $n$, dimension $k$, minimum distance $d$ and covering
radius $R$ is denoted by $[n,k,d]_qR$.

A coset of the code $C$ defined by the vector $x \in F_q^n$ is the
set $x+C=\{x+c \ \vert \ c \in C\}$. A coset leader of $x+C$ is a
vector in $x+C$ of smallest weight. When the code is linear its
covering radius is equal to the weight of the heaviest coset
leader. The covering radius of a linear code can also be defined
in terms of the parity check matrix.

\begin{theorem} \cite{CHLL}
Let $C$ be a $[n,k]$ code with parity check matrix $H$. The
covering radius of $C$ is the smallest integer $R$ such that every
$q$-ary $(n-k)$-tuple can be written as a linear combination of at
most $R$ columns of $H$.
\end{theorem}

The special case are codes for which $R(C)=e(C)$ and such codes
are called {\it perfect} codes. The problem of finding all perfect
codes was begun by Golay in 1949 and completed in 1973 by
Tiet\"av\"ainen \cite{T} and independently by Zinov'ev and
Leont'ev \cite{ZL}. The only perfect codes are: $[n,n,1]_q0$ codes
for each $n\geq 1$; $[2s+1,1,2s+1]_2s$ repetition codes for each
$s\geq 1$; code of length $n$ containing only one codeword;
$q$-ary codes with the parameters of Hamming codes; the
$[23,12,7]_23$ binary Golay code; the $[11,6,5]_32$ ternary Golay
code.

The next step in this direction is to consider codes for which
packing and covering radii differ by 1, i.e.\ {\it quasi-perfect}
codes. A code is called quasi-perfect (QP) if its packing radius
is $e$ and its covering radius is $e+1$, for some nonnegative
integer $e$. Clearly, the minimum distance of such a code is
$2e+1$ or $2e+2$. Then a natural question is which codes are
quasi-perfect? It is clear that any code with covering radius 1
and minimum distance 1 or 2 is quasi-perfect. Therefore,
quasi-perfect codes with covering radius 1 are not interesting and
we will focus on the investigation of quasi-perfect codes with
covering radius greater than 1.

\section {Known results about quasi-perfect codes with covering
radius greater than 1}

Quasi-perfect codes with covering radius 2 and 3 were extensively
studied and many infinite families of binary, ternary and
quaternary QP codes are known. In particular, codes with
parameters $[n,k,d]_q2$ $d=3,4$ are QP. These codes are connected
with 1-saturating sets in projective spaces $PG(n-k-1,q)$ and a
lot of infinite families of such codes are described in the
literature (see \cite{GS} - \cite{GP}). The following theorem
leads to a chain of QP codes.

\begin{theorem}
Assume that an $[n,k,d]_q2$ QP code with ${\displaystyle n\leq
{{q^{n-k}-1}\over {q-1}} - 2}$ and $3\leq d\leq 4$ exists. Then
an$[n+1,k+1,3]_q2$ QP code exists.
\end{theorem}

{\it Proof.} Let we add a column to a parity check matrix of a
$[n,k,d]_q2$ code to obtain a new $[n+1,k+1,d_{new}]_qR_{new}$
code. According to the Theorem 1 $R_{new}\leq 2$. As the new
length ${\displaystyle n+1\leq {{q^{n-k}-1}\over {q-1}} - 1}$, it
is impossible that $R_{new}=1$. Also for any new column it holds
that $d_{new}\leq 3$. If the new column is not obtained by
multiplication of an old column by an element of the field $GF(q)$
then $d_{new}=3$. The right choice of the new column is possible
as ${\displaystyle n\leq {{q^{n-k}-1}\over {q-1}} - 2}$.

$[n,n-4,5]_q3$ QP codes correspond to complete arcs in the
projective space $PG(3,q)$ and are investigated in \cite{DMP1},
\cite{DFMP}. Survey and new results for the binary QP codes can be
found in \cite{EM}. It appears that there is a great variety of QP
codes of minimum distance up to 5.

Considerably less is known for $q$-ary QP codes with $q>2$. One
infinite family of ternary codes is known due to Gashkov and
Sidel'nikov \cite{GS}. The family members are
$[(3^s+1)/2,(3^s+1)/2-2s,5]_33$ codes. Quasi-perfectness of two
families of quaternary codes, namely
$[(4^s-1)/3,(4^s-1)/3-2s,5]_43$ and
$[(2^{2s+1}+1)/3,(2^{2s+1}+1)/3-2s-1,5]_43$, presented in \cite{G}
and \cite{DZ} was shown by Dodunekov \cite{Dod,Dod1}.

The first computer searches to find new quasi-perfect codes were
by Wagner in 1966~\cite{W}. He proposed a tree-search program
which uses the properties of parity-check matrices of binary
linear quasi-perfect codes to find such codes. Fixing the number
of check digits and the number of errors to be corrected, the
program finds one quasi-perfect code for each block length if such
a code exists. Using this program 27 new binary linear QP codes
were found \cite{W,W1}. The codes have lengths between 19 and 55
and all have covering radius 3. Later Simonis \cite{Sim} proved
that one of Wagner's codes, namely $[23,14,5]$, is unique.

Baicheva, Dodunekov and K\"otter \cite{BDK} investigated the
weight structure and error-correcting performance of the ternary
$[13,7,5]$ quadratic-residue code and showed that the covering
radius of the code is equal to three, i.e.\ it is a quasi-perfect
code. Recently Danev and Dodunekov \cite{DD} proved that this code
is the first member of a family of ternary QP codes with
parameters $[(3^s-1)/2,(3^s-1)/2-2s,5]_3$ for all odd $s\geq 3$.

All these results lead to the following question: How restrictive
is quasi-perfectness, i.e.\ are there inequivalent quasi-perfect
codes? In our work we classify all binary of dimension up to 9 and
ternary of dimension up to 6 linear QP codes as well as give some
partial classifications for dimensions up to 14 and 13
respectively. It turned out that there are many cases where more
than one QP code for fixed length and dimension exists. In this
way we answer the above question.

\section{Classification of binary and ternary quasi-perfect
linear codes}

The approach used in this work is based on the classification of
codes with given parameters. First we fix the dimension of the
code and determine the possible lengths and minimum distances of
the codes which could be quasi-perfect. Then we classify all such
codes and finally compute their covering radii. In this way we
determine all quasi-perfect codes with the fixed parameters. In
order to determine the parameters of possible candidates for
quasi-perfect codes with covering radius $e+1$ we take into
account that the minimum distance of these codes can only be
$2e+1$ or $2e+2$. Brouwer's tables of bounds on the size of linear
codes \cite{Bro} and the tables for the least covering radius of
binary \cite{CHLL} and ternary linear codes \cite{BV} are used to
find the possible lengths of QP codes when the dimension is fixed.
Once the parameters (length, dimension and minimum distance) are
determined all codes with these parameters are classified up to
equivalence using the approach of \cite{Bou1}.

In the classification of the codes two main approaches were used.
The first one is based on puncturing, the second one on
shortening. While in general the dimension of the code is
unchanged by puncturing, this is not true if all non-zero
positions of a codeword are deleted. Let $G$ be a generator matrix
of a linear $[n,k,d]_q$ code $C$. Then the residual code
$\mbox{Res}(C,{\bf c})$ of $C$ with respect to a codeword $c$ is
the code generated by the restriction of $G$ to the columns where
$c$ has a zero entry. A lower bound on the minimum distance of the
residual code is given by

\begin{lemma}  \cite{D}
Suppose $C$ is an $[n,k,d]_q$ code and suppose ${\bf c} \in C$ has
weight w, where $d > w(q-1)/q$. Then $\mbox{Res}(C,{\bf c})$ is an
$[n-w,k-1,d']_q$ code with $d' \geq  d - w + \lceil{w/q}\rceil$.
\end{lemma}

Inverting this operation, we search for an $[n,k,d]_q$ code on the
basis of an $[n-w,k-1,d']_q$ code (its residual code with respect
to a codeword of weight $w$) or an $[n-i,k,d']_q$ code (punctured
on $i$ coordinates code). We can apply the same operation for the
residual or punctured code respectively. This procedure is
repeated until we obtain as a start code a code with small
parameters such that all codes having these parameters can easily
be classified. For example, starting from the $[3,2,2]_2$ code, we
obtain all $[8,3,5]_2$ codes, take only the nonequivalent of them
and then obtain all nonequivalent $[28,4,20]_2$ codes.

The second approach increases both the length and the dimension of
the code, i.e.\ we construct $[n,k,d]_q$ codes extending
\mbox{$[n{-}i,k{-}i,d]_q$} or $[n-i-1,k-i,d]_q$ codes. The
following result shows when the latter type of code can be used
\cite[p.\ 592]{MS}.

\begin{lemma}
Let $C$ be an $[n,k,d]_q$ code. If there exists a codeword ${\bf
c} \in C^{\perp}$ with $wt({\bf c})=i$, then there is an
$[n-i,k-i+1,d]_q$ code.
\end{lemma}

If $G$ is a generator matrix for an $[n-i,k-i,d]_q$ or an
$[n-i-1,k-i,d]_q$ code, we extend it (in all possible ways) to
\begin{equation}
\left(
\begin{array}{c|c}
* & {\bf I}_{i} \\ \hline
{\bf G} & {\bf 0}
\end{array}
\right) \ \ \rm{or} \ \ \left(
\begin{array}{c|c}
* & {\bf 1} \ {\bf I}_{i} \\ \hline
{\bf G} & {\bf 0}
\end{array}
\right),
\end{equation}
respectively, where ${\bf I}_{i}$ is the $i\times i$ identity
matrix, {\bf 1} is an all-1 column vector, and the starred
submatrix is to be determined. If we let the matrix $G$ be in
systematic form, we can fix $k$ more columns to get
\begin{equation}
\left(
\begin{array}{c|c|c}
* & {\bf 0} & {\bf I}_{i} \\ \hline
{\bf G}_{1} & {\bf I}_{k} & {\bf 0}
\end{array}
\right) \ \ \rm{or} \ \ \left(
\begin{array}{c|c|c}
* & {\bf 0} & {\bf 1} \ {\bf I}_{i} \\ \hline
{\bf G}_{1} & {\bf I}_{k} & {\bf 0}
\end{array}
\right).
\end{equation}

Again we can apply recursively the same approach to obtain $G_1$
while on the bottom of this hierarchy of extensions is the trivial
$[k,k,1]_q$ code.

 In our investigation we are interested in codes with covering radius
greater than 1, therefore they have minimum distance at least 3.
Thus their dual codes are projective codes. To classify the binary
codes with codimension $n-k$ up to 6 we use the results from
\cite{Bou} where all binary projective codes with dimensions up to
6 are classified. Then among the codes from \cite{Bou} we consider
only those having the necessary minimum distance of the dual code.
For example, to classify all $[8,2,5]$ codes we consider 14
$[8,6]$ projective codes. The dual code of only one of them has
minimum distance 5 and thus we have only one $[8,2,5]$ code. In
the same way using results from \cite{BB} where all ternary
projective codes of dimension 4 are classified we classify the
ternary linear codes of codimension 4 which could be
quasi-perfect.

After the classification was completed we proceed with the
determination of the covering radii of the codes. We recall that
in the case of linear codes this is equivalent to the
determination of the heaviest coset leader. To do this, we use the
fact that if the code is in a systematic form, a representative of
each coset can be found by generating all words of the form
$(\underbrace {0,\dots ,0}_{k},a)$, $a\in F_q^{n-k}$. Taking into
account that each vector of weight less than or equal to $e$ is a
unique coset leader, we test only words of the above form and
weight greater than $e$. Therefore we have to test at most
$\sum_{i=e}^{n-k}{{n-k}\choose{i+1}}(q-1)^{i+1}$ words because if
we obtain a coset leader of weight greater than $e+1$ we stop the
check.

{\it Remark.} Let us denote by $\alpha_i$ for $i=0,1,\dots ,n$ the
number of coset leaders of weight $i$. The set of the coset
leaders of each QP code is known. All vectors of weights less than
or equal to $e$ are coset leaders and thus $\alpha_i={n\choose
i}(q-1)^i$ for $i=0,\dots ,e$. Then for $\alpha_{e+1}$ we get
$\alpha_{e+1}=q^k-\sum_{i=0}^e\alpha_i$.

\section{Results}
By the approach described in the previous section all binary and
ternary quasi-perfect codes of dimensions up to 9 and 6
correspondingly are determined, as well as some partial results
for binary codes of dimensions up to 14 and ternary codes of
dimensions up to 13 are obtained. The results are summarized in
Table I.

\begin{table}
\caption{Binary and Ternary Quasi-perfect Codes} \centering
\begin{tabular}{|l|l|l||l|l|l|} \hline
\multicolumn{6}{|c|}{Binary quasi-perfect codes}\\ \hline
Code&All&QP&Code&All&QP\\
\hline
$[5,2,3]$ &  1 & 1&$[14,9,3]$ & 126 &113 \\
$[6,3,3]$ &  1 & 1&$[15,9,3]$& 11464& 380\\
$[8,2,5]$* &  1 & 1&$[17,9,5]$& 1  & 1\\
$[7,3,3]$ &  3 &2&$[14,10,3]$&  1& 1\\
$[8,4,4]$*&  1 & 1&$[15,10,3]$& 142&131 \\
$[8,4,3]$ &  4 & 4&$[16,10,3]$& 28900& 2296\\
$[9,4,4]$* &  4 & 1&$[19,10,5]$& 31237 & 13\\
$[9,4,3]$ &  19 & 1&$[16,11,4]$*& 1 & 1\\
$[11,4,5]$&  1 & 1&$[16,11,3]$& 143 & 143\\
$[9,5,3]$ &  5 & 5 &$[17,11,4]$*& 39 & 5\\
$[10,5,4]$* & 4 & 1&$[17,11,3]$&70416& 12221\\
$[10,5,3]$& 37  & 12&$[20,11,5]$&13924 & 565\\
$[10,6,3]$& 4 & 4&$[17,12,3]$& 129 & 129\\
$[11,6,3]$& 58& 25&$[18,12,4]$*& 33 & 1\\
$[14,6,5]$& 11 & 1&$[21,12,5]$&2373 & 666\\
$[11,7,3]$& 3  & 3&$[22,12,6]$& 128  & 1\\
$[12,7,3]$& 84 & 55&$[24,12,8]$& 1  & 1\\
$[13,7,4]$*& 45  & 1&$[24,12,7]$&  11& 11\\
$[13,7,3]$& 1660  & 7&$[25,12,8]$&  7 & 2\\
$[15,7,5]$&  6 & 4&$[18,13,3]$& 113& 113 \\
$[12,8,3]$&   2& 2&$[19,13,3]$& 366064& 185208\\
$[13,8,3]$& 109 & 88&$[22,13,5]$&128 & 120\\
$[14,8,3]$& 4419 & 65&$[19,14,3]$&  91 & 91\\
$[13,9,3]$& 1  & 1&$[20,14,4]$*&  24 & 1\\

 \hline \hline
 \multicolumn{6}{|c|}{Ternary quasi-perfect codes}
 \\ \hline
Code&All&QP&Code&All&QP\\
\hline
$[5,2,3]$* &  2 &2&$[11,7,3]$&  339  & 319\\
$[6,3,3]$ &1 & 1&$[12,7,3]$& 60910 & 1\\
$[7,4,3]$ &  4 &4& $[13,7,5]$&  6  & 5\\
$[8,4,4]$* & 3  &2&$[11,8,3]$&  1  & 1\\
$[8,4,3]$ & 37  &5&$[12,8,3]$& 805 & 753\\
$[8,5,3]$ & 3  &3&$[14,8,5]$ &1 &1\\
$[9,5,3]$ & 87 &23&$[12,9,3]$ &1 &1\\
$[9,6,3]$ & 3  &3&$[13,9,3]$ &1504 &1479\\
$[10,6,4]$*&  1 &1&$[14,10,3]$ &2695 &2659\\
$[10,6,3]$& 195&102&$[15,11,3]$ &4304 &4304\\
$[12,6,6]$& 1   &1&$[16,12,3]$ &6472 &6472\\
$[12,6,5]$& 36  &18&$[17,13,3]$ &8846 &8846\\
$[10,7,3]$& 2 & 2 & & &\\
\hline
\end{tabular}
\end{table}

Some of the codes from the table are not new and have already been
constructed in previous works. We will note that QP codes with
minimum distances 3 or 4 and covering radius 2 are connected with
1-saturating sets in projective spaces $PG(n-k-1,q)$ in the
following way: the points of a 1-saturating $n$-set can be
considered as $n-k$-dimensional columns of a parity-check matrix
of an $[n,k]_q2$ code. Also QP codes with minimum distance 4 are
complete caps in $PG(n-k-1,q)$. Constructions of minimal
1-saturating sets and complete caps in binary projective spaces
$PG(k-1,2)$ are described in \cite{FMMP}, \cite{Dav}, \cite{DMP} -
\cite{DMP2}, \cite{KL}. Codes obtained in these works are marked
with a *. Some of the marked codes are also obtained in \cite{KL}
where recursive constructions of complete caps in $PG(n-k-1,2)$
are given. Existing of codes with parameters $[10,5,3]_22$,
$[14,6,5]_23$ and $[13,7,3]_22$ is shown in \cite{GS}. The
$[17,9,5]_23$ code is the first representative from the infinite
family of $[2^{2s}+1,2^{2s}+1-4,5]_2$, $s\geq 2$ Zetterberg's
codes which are proved to be quasi-perfect by Dodunekov
\cite{Dod}. $[19,10,5]_23$, $[20,11,5]_23$, $[23,14,5]_23$ and
$[24,14,6]_23$ are among the QP codes obtained by a computer
search by Wagner. He obtained only one representative for each of
the parameters. Our classification shows that QP codes with the
first two parameters are not unique. There are additionally 12
$[19,10,5]_2$ and 564 $[20,11,5]_2$ quasi-perfect codes.
$[24,12,8]_2$ is the well known extended Golay code which is also
known to be a quasi-perfect one.

For a completeness of the classification results about QP codes,
we will note some not classified in this work such codes. In
\cite{DFMP} the unique $[6,1,5]_23$ and in \cite{DMP1} the unique
$[5,1,5]_33$ codes are presented. As complete caps in $PG(4,3)$
the $[16,11,4]_32$, $[17,12,4]_32$ and $[18,13,4]_32$ codes in
\cite{FMMP} and in $PG(6,2)$ the $[21,14,4]_22$ code in \cite{GDT}
are obtained. In \cite{KL} it is shown that there are 5
nonequivalent $[21,14,4]_22$ codes. Also applying Theorem 2 to
codes from the Table the following chains of QP codes' parameters
can be obtained.\\
$[5,2,3]_22\rightarrow [6,3,3]_22;$\\
$[8,4,4]_22\rightarrow \dots \rightarrow [14,10,3]_22;$\\
$[9,4,4]_22\rightarrow \dots \rightarrow [30,25,3]_22;$\\
$[13,7,4]_22\rightarrow \dots \rightarrow {\bf
[18,12,3]_22}\rightarrow [19,13,3]_22\rightarrow {\bf
[20,14,3]_22}\rightarrow \dots \rightarrow
{\bf[62,56,3]_22};$\\
$[5,2,3]_32\rightarrow [12,3,3]_32;$\\
$[8,4,4]_32\rightarrow \dots \rightarrow [17,13,3]_32\rightarrow
{\bf [18,14,3]_32}\rightarrow \dots \rightarrow {\bf
[40,36,3]_32}$;\\
$[12,7,3]_32\rightarrow {\bf [13,8,3]_32}\rightarrow \dots \rightarrow {\bf [121,116,3]_32}$.\\
 Codes not classified in this work are boldfaced.

Until this work the only known examples of QP codes with minimum
distance greater than 5 were binary repetition codes, the
$[24,12,8]_24$ extended Golay code, the $[22,12,6]_23$ punctured
Golay code, $[7,1,7]_34$ and $[8,1,7]_24$ codes classified in
\cite{DFMP}. We provide examples of more such codes and in this
way answer the first open question from the recent paper of Etzion
and Mounits \cite{EM} where to find new or to prove the
nonexistence of QP codes with $d>5$ is suggested. The most
interesting are $[24,12,7]_24$ and $[25,12,8]_24$ codes  which are
the first examples of quasi-perfect codes with $R=4$ except the
$[24,12,8]_24$ extended Golay and the $[8,1,8]_24$ repetition
codes. The generator matrices of these codes are given in the
Appendix. The codes are in a systematic form with generator matrix
$G=[I_k\vert A]$ and the identity matrix $I_k$ is omitted in order
to save space.

\section{Conclusions}

In this work classification results about binary and ternary
linear quasi-perfect codes of small dimensions are obtained. More
precisely, all binary QP codes of dimensions up to 9 and ternary
QP codes of dimensions up to 6 are classified as well as some
partial classifications about QP codes of dimensions up to 14 are
got. The results show that for each dimension there are only few
possible lengths for which quasi-perfect codes exist. For some
parameters hundreds and thousands of nonequivalent QP codes are
found which means that quasi-perfectness is not so restrictive
characteristic of the code. QP codes of minimum distance greater
than 5 are obtained and therefore it could be expected that at
greater dimensions QP codes with bigger covering radii exist. Thus
it will be an interesting research problem to answer the following
questions:

$\bullet$ Are there quasi-perfect codes with minimum distance
greater than 8 except the binary repetition code?

$\bullet$ Is there an upper bound about minimum distance of a QP
code?

At the end we will conclude with the observation that the
classification of all parameters of QP codes would be much more
difficult than the similar one for perfect codes.

 \section {Acknowledgement} The authors wish to express her appreciation to the anonymous
reviewers whose comments and suggestions greatly improved the
paper. Theorem 2 from section II is due to one of the reviewers.

\section {Appendix. Generator matrices of binary quasi-perfect codes
with minimum distance 7 and 8}

\centerline{A. $[24,12,7]_24$ QP codes} {\small $A_1=\left(
\begin{array}{c}
010010110101\\ 101000111001\\
111010010010\\ 000110011011\\
010101010110\\ 101100001110\\
111111011101\\ 111100110100\\
101101010011\\ 010101111001\\
011000011111\\ 000011111110 \end{array} \right )$ $A_2=\left(
\begin{array}{c}
110001001011\\ 001011000111\\
111010010010\\ 100101100101\\
010101010110\\ 101100001110\\
011100100011\\ 111100110100\\
001110101101\\ 110110000111\\
111011100001\\ 000011111110\end{array} \right )$

$A_3=\left( \begin{array}{c}
010101001011\\ 001100111001\\
011110010010\\ 100010011011\\
110110101000\\ 101100001110\\
111111011101\\ 111011001010\\
101010101101\\ 110110000111\\
011000011111\\ 000011111110\end{array} \right )$ $A_4=\left(
 \begin{array}{c}
100010111101\\ 001011010011\\
101110001110\\ 110011101010\\
110111010100\\ 111101111111\\
111110100001\\ 011011001101\\
111010110110\\ 110110011011\\
001111111000\\ 000111100111\end{array} \right ) $

$A_5=\left(  \begin{array}{c}
000010111101\\ 011100110100\\
011001101001\\ 110100001101\\
110000110011\\ 101101111111\\
101110100001\\ 111100101010\\
111101010001\\ 100110011011\\
001111111000\\ 010111100111\end{array} \right )$ $A_6=\left(
 \begin{array}{c}
000011011110\\ 101010101011\\
101110010101\\ 010011110001\\
010110101100\\ 111101111111\\
111111000010\\ 111011001101\\
111010110110\\ 010110011011\\
001111111000\\ 000111100111\end{array} \right )$

$A_7=\left(  \begin{array}{c}
110001010101\\ 001011011001\\
111010010010\\ 000110010111\\
110110100100\\ 001111100010\\
111111001111\\ 111100101010\\
001110101101\\ 110110011001\\
111011100001\\ 000011111110\end{array} \right ) $ $A_8=\left(
 \begin{array}{c}
000011011011\\ 001010101101\\
001110010110\\ 100011100110\\
100110110001\\ 101101111111\\
101111001000\\ 101011010101\\
101010111010\\ 100110001111\\
001111100011\\ 010111111100\end{array} \right ) $

$A_9=\left(  \begin{array}{c}
010101011101\\ 011011101100\\
001110001111\\ 100001101011\\
110111000110\\ 111100100001\\
101010110010\\ 101101010111\\
111110011100\\ 110011110101\\
011000111011\\ 000111111010\end{array} \right )A_{10}=\left(
 \begin{array}{c}
110001101110\\ 111111011111\\
101010111100\\ 010110101101\\
100100110011\\ 001011100111\\
011101110100\\ 111110100010\\
101101101001\\ 010111000110\\
011000111011\\ 000111111010\end{array} \right )$

$A_{11}=\left( \begin{array}{c}
000111100110\\ 101101010101\\
101100101011\\ 110000110111\\
110011011100\\ 111011100001\\
111110010010\\ 011010101110\\
011001011011\\ 010101101101\\
001111111000\\ 000110011111\end{array} \right )$

\vspace{4cm}

\centerline{ B. $[25,12,8]_24$ QP codes}

$A_1=\left( \begin{array}{c}1101101100100\\
1101000111001\\
1110100001101\\ 1110110110000\\
1011001110010\\ 1011010101100\\
1000100111110\\ 0111100101010\\
0111101010001\\ 0100110011011\\
0001111111000\\ 0010111100111\end{array}
\right ) $ $A_2=\left( \begin{array}{c}1101101110000\\
1101000100111\\
1100011101001\\ 1100001011110\\
1001110010110\\ 1001101001101\\
1000100111011\\ 0101011010101\\
0101010111010\\ 0100110001111\\
0001111100011\\ 0010111111100\end{array} \right)$ }

\end{document}